# Icosahedral Polyhedra from $D_6$ lattice and Danzer's *ABCK* tiling

Abeer Al-Siyabi, [a] Nazife Ozdes Koca, [a*] and Mehmet Koca[b]

[a]Department of Physics, College of Science, Sultan Qaboos University, P.O. Box 36, Al-Khoud, 123 Muscat, Sultanate of Oman, [b]Department of Physics, Cukurova University, Adana, Turkey, retired, *Correspondence e-mail: nazife@squ.edu.om

## ABSTRACT

It is well known that the point group of the root lattice $D_6$ admits the icosahedral group as a maximal subgroup. The generators of the icosahedral group $H_3$, its roots and weights are determined in terms of those of $D_6$. Platonic and Archimedean solids possessing icosahedral symmetry have been obtained by projections of the sets of lattice vectors of $D_6$ determined by a pair of integers $(m_1, m_2)$ in most cases, either both even or both odd. Vertices of the Danzer's *ABCK* tetrahedra are determined as the fundamental weights of $H_3$ and it is shown that the inflation of the tiles can be obtained as projections of the lattice vectors characterized by the pair of integers which are linear combinations of the integers $(m_1, m_2)$ with coefficients from Fibonacci sequence. Tiling procedure both for the *ABCK* tetrahedral and the $<ABCK>$ octahedral tilings in 3D space with icosahedral symmetry $H_3$ and those related transformations in 6D space with $D_6$ symmetry are specified by determining the rotations and translations in 3D and the corresponding group elements in $D_6$. The tetrahedron *K* constitutes the fundamental region of the icosahedral group and generates the rhombic triacontahedron upon the group action. Properties of "*K*-polyhedron", "*B*-polyhedron" and "*C*-polyhedron" generated by the icosahedral group have been discussed.

# 1. Introduction

Quasicrystallography as an emerging science attracts the interests of many scientists varying from the fields of material science, chemistry and physics. For a review see for instance the references (Di Vincenzo & Steinhardt, 1991; Janot, 1993; Senechal, 1995). It is mathematically intriguing as it requires the aperiodic tiling of the space by some prototiles. There have been several approaches to describe the aperiodicity of the quasicrystallographic space such as the set theoretic techniques, cut-and-project scheme of the higher dimensional lattices and the intuitive approaches such as the Penrose-like tilings of the space. For a review of these techniques, we propose the reference (Baake & Grimm, 2013).

There have been two major approaches for the aperiodic tiling of the 3D space with local icosahedral symmetry. One of them is the Socolar-Steinhardt tiles (Socolar & Steinhardt, 1986) consisting of acute rhombohedron with golden rhombic faces, Bilinski rhombic dodecahedron, rhombic icosahedron and rhombic triacontahedron, the latter three are constructed with two Ammann tiles of acute and obtuse rhombohedra. Later it was proved that (Danzer, Papadopolos & Talis, 1993), (Roth, 1993) they can be constructed by the Danzer's *ABCK* tetrahedral tiles (Danzer, 1989). Katz (Katz, 1989) and recently Hann-Socolar-Steinhardt (Hann, Socolar & Steinhardt, 2018) proposed a model of tiling scheme with decorated Ammann tiles. A detailed account of the Danzer *ABCK* tetrahedral tilings can be found in "math.uni-bielefeld.de/icosahedral tilings in $\mathbb{R}^3$: the *ABCK* tilings" and in page 231 of the reference (Baake & Grimm, 2013) where the substitution matrix, its eigenvalues and the corresponding eigenvectors are studied. The right and left eigenvectors of the substitution matrix corresponding to the Perron-Frobenius eigenvalue are well known and will not be repeated here.

Ammann rhombohedral and Danzer *ABCK* tetrahedral tilings are intimately related with the projection of six-dimensional cubic lattice and the root and weight lattices of $D_6$, the point symmetry group of which is of order $2^5 6!$. See for a review the paper "Modelling of quasicrystals" by Kramer (Kramer, 1993) and references therein. Similar work has also been carried out in the reference (Koca, Koca & Koc, 2015). Kramer and Andrle (Kramer & Andrle, 2004) have investigated Danzer tiles from the wavelet point of view and their relations with the lattice $D_6$.

In what follows we point out that the icosahedral symmetry requires a subset of the root lattice $D_6$ characterized by a pair of integers $(m_1, m_2)$ with $m_1 + m_2 = even$ which are the coefficients of the orthogonal set of vectors $l_i, (i = 1, 2, ..., 6)$. Our approach is different than the cut and project scheme of lattice $D_6$ as will be seen in the sequel. The paper consists of two major parts; first part deals with the determination of Platonic and Archimedean icosahedral polyhedra by projection of the fundamental weights of the root lattice $D_6$ into 3D space and the second part employs the technique to determine the images of the Danzer tiles in $D_6$. Inflation of the Danzer tiles are related to a redefinition of the pair of integers $(m_1, m_2)$ by the Fibonacci sequence. Embeddings of basic tiles in the inflated ones require translations and rotations in 3D space where the corresponding transformations in 6D space can be easily determined by the technique we have introduced. This technique which restricts the lattice $D_6$ to its subset has not been discussed elsewhere.

The paper is organized as follows. In Sec. 2, we introduce the root lattice $D_6$, its icosahedral subgroup, decomposition of its weights in terms of the weights of the icosahedral group $H_3$ leading to the Archimedean polyhedra projected from the $D_6$ lattice. It turns out that the lattice vectors to be projected are determined by a pair of integers $(m_1, m_2)$. In Sec. 3 we introduce the *ABCK* tetrahedral Danzer tiles in terms of the fundamental weights of the icosahedral group $H_3$. Tiling by inflation with $\tau^n$ where $\tau = \frac{1+\sqrt{5}}{2}$ and $n \in \mathbb{Z}$ is studied in $H_3$ space by prescribing the appropriate rotation and translation operators. The corresponding group elements of $D_6$ are determined noting that the pair of integers $(m_1, m_2)$ can be expressed as the linear combinations of similar integers with coefficients from Fibonacci sequence. Sec. 4 is devoted for conclusive remarks.

## 2. Projection of $D_6$ lattice under the icosahedral group $H_3$ and the Archimedian polyhedra

We will use the Coxeter diagrams of $D_6$ and $H_3$ to introduce the basic concepts of the root systems, weights and the projection technique. A vector of the $D_6$ lattice can be written as a linear combination of the simple roots with integer coefficients:

$$\lambda = \sum_{i=1}^{6} n_i \alpha_i = \sum_{i=1}^{6} m_i l_i, \quad \sum_{i=1}^{6} m_i = even, \quad n_i, m_i \in \mathbb{Z}. \tag{1}$$

Here $\alpha_i$ are the simple roots of $D_6$ defined in terms of the orthonormal set of vectors as $\alpha_i = l_i - l_{i+1}, i = 1, \ldots, 5$ and $\alpha_6 = l_5 + l_6$ and the reflection generators of $D_6$ act as $r_i: l_i \leftrightarrow l_{i+1}$ and $r_6: l_5 \leftrightarrow -l_6$. The generators of $H_3$ can be defined as (Koca, Koca & Koc, 2015) $R_1 = r_1 r_5, R_2 = r_2 r_4$ and $R_3 = r_3 r_6$ where the Coxeter element, for example, can be taken as $R = R_1 R_2 R_3$. The weights of $D_6$ are given by

$$\omega_1 = l_1, \ \omega_2 = l_1 + l_2, \ \omega_3 = l_1 + l_2 + l_3, \ \omega_4 = l_1 + l_2 + l_3 + l_4,$$
$$\omega_5 = \frac{1}{2}(l_1 + l_2 + l_3 + l_4 + l_5 - l_6), \ \omega_6 = \frac{1}{2}(l_1 + l_2 + l_3 + l_4 + l_5 + l_6). \tag{2}$$

The Voronoi cell of the root lattice $D_6$ is the dual polytope of the root polytope of $\omega_2$ (Koca et. al, 2018) and determined as the union of the orbits of the weights $\omega_1$, $\omega_5$ and $\omega_6$ which correspond to the holes of the root lattice (Conway & Sloane, 1999). If the roots and weights of $H_3$ are defined in two complementary 3D spaces $E_\parallel$ and $E_\perp$ as $\beta_i(\acute{\beta}_i)$ and $v_i(\acute{v}_i), (i = 1, 2, 3)$ respectively then they can be expressed in terms of the roots and weights of $D_6$ as

$$\beta_1 = \frac{1}{\sqrt{2+\tau}}(\alpha_1 + \tau\alpha_5), \ v_1 = \frac{1}{\sqrt{2+\tau}}(\omega_1 + \tau\omega_5),$$
$$\beta_2 = \frac{1}{\sqrt{2+\tau}}(\alpha_2 + \tau\alpha_4), \ v_2 = \frac{1}{\sqrt{2+\tau}}(\omega_2 + \tau\omega_4), \tag{3}$$
$$\beta_3 = \frac{1}{\sqrt{2+\tau}}(\alpha_6 + \tau\alpha_3), \ v_3 = \frac{1}{\sqrt{2+\tau}}(\omega_6 + \tau\omega_3).$$

For the complementary 3D space replace $\beta_i$ by $\acute{\beta}_i$, $v_i$ by $\acute{v}_i$ in (3) and $\tau = \frac{1+\sqrt{5}}{2}$ by its algebraic conjugate $\sigma = \frac{1-\sqrt{5}}{2} = -\tau^{-1}$. Consequently, Cartan matrix of $D_6$ (Gram matrix in lattice terminology) and its inverse block diagonalize as

$$C = \begin{bmatrix} 2 & -1 & 0 & 0 & 0 & 0 \\ -1 & 2 & -1 & 0 & 0 & 0 \\ 0 & -1 & 2 & -1 & 0 & 0 \\ 0 & 0 & -1 & 2 & -1 & -1 \\ 0 & 0 & 0 & -1 & 2 & 0 \\ 0 & 0 & 0 & -1 & 0 & 2 \end{bmatrix} \rightarrow \begin{bmatrix} 2 & -1 & 0 & 0 & 0 & 0 \\ -1 & 2 & -\tau & 0 & 0 & 0 \\ 0 & -\tau & 2 & 0 & 0 & 0 \\ 0 & 0 & 0 & 2 & -1 & 0 \\ 0 & 0 & 0 & -1 & 2 & -\sigma \\ 0 & 0 & 0 & 0 & -\sigma & 2 \end{bmatrix}$$
$$\tag{4}$$

$$C^{-1} = \frac{1}{2}\begin{bmatrix} 2 & 2 & 2 & 2 & 1 & 1 \\ 2 & 4 & 4 & 4 & 2 & 2 \\ 2 & 4 & 6 & 6 & 3 & 3 \\ 2 & 4 & 6 & 8 & 4 & 4 \\ 1 & 2 & 3 & 4 & 3 & 2 \\ 1 & 2 & 3 & 4 & 2 & 3 \end{bmatrix} \rightarrow \frac{1}{2}\begin{bmatrix} 2+\tau & 2\tau^2 & \tau^3 & 0 & 0 & 0 \\ 2\tau^2 & 4\tau^2 & 2\tau^3 & 0 & 0 & 0 \\ \tau^3 & 2\tau^3 & 3\tau^2 & 0 & 0 & 0 \\ 0 & 0 & 0 & 2+\sigma & 2\sigma^2 & \sigma^3 \\ 0 & 0 & 0 & 2\sigma^2 & 4\sigma^2 & 2\sigma^3 \\ 0 & 0 & 0 & \sigma^3 & 2\sigma^3 & 3\sigma^2 \end{bmatrix}.$$

Fig.1 shows how one can illustrate the decomposition of the 6D space as the direct sum of two complementary 3D spaces where the corresponding nodes $\alpha_i, (i = 1, 2, \ldots, 6)$, $\beta_i, (i = 1, 2, 3)$ represent the corresponding simple roots and $\acute{\beta}_i$ stands for $\beta_i$ in the complementary space.

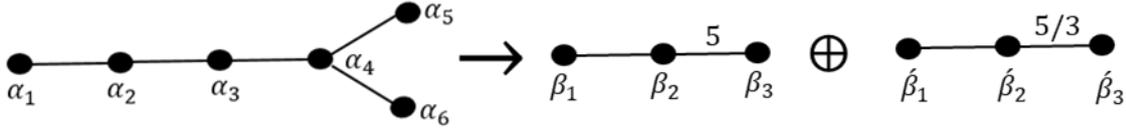

**Figure 1**
Coxeter-Dynkin diagrams of $D_6$ and $H_3$ illustrating the symbolic projection.

For a choice of the simple roots of $H_3$ as $\beta_1 = (\sqrt{2}, 0, 0), \beta_2 = -\frac{1}{\sqrt{2}}(1, \sigma, \tau), \beta_3 = (0, 0, \sqrt{2})$ and the similar expressions for the roots in the complementary space we can express the components of the set of vectors $l_i$ $i = 1, 2, \ldots, 6$) as

$$\begin{bmatrix} l_1 \\ l_2 \\ l_3 \\ l_4 \\ l_5 \\ l_6 \end{bmatrix} = \frac{1}{\sqrt{2(2+\tau)}} \begin{bmatrix} 1 & \tau & 0 & \tau & -1 & 0 \\ -1 & \tau & 0 & -\tau & -1 & 0 \\ 0 & 1 & \tau & 0 & \tau & -1 \\ 0 & 1 & -\tau & 0 & \tau & 1 \\ \tau & 0 & 1 & -1 & 0 & \tau \\ -\tau & 0 & 1 & 1 & 0 & \tau \end{bmatrix}. \tag{5}$$

First three and last three components project the vectors $l_i$ into $E_\parallel$ and $E_\perp$ spaces respectively. On the other hand the weights are represented by the vectors $v_1 = \frac{1}{\sqrt{2}}(1, \tau, 0), v_2 = \sqrt{2}(0, \tau, 0), v_3 = \frac{1}{\sqrt{2}}(0, \tau^2, 1)$. The generators of $H_3$ in the space $E_\parallel$ read

$$R_1 = \begin{bmatrix} -1 & 0 & 0 \\ 0 & 1 & 0 \\ 0 & 0 & 1 \end{bmatrix}, \quad R_2 = \frac{1}{2}\begin{bmatrix} 1 & -\sigma & -\tau \\ -\sigma & \tau & 1 \\ -\tau & 1 & \sigma \end{bmatrix}, \quad R_3 = \begin{bmatrix} 1 & 0 & 0 \\ 0 & 1 & 0 \\ 0 & 0 & -1 \end{bmatrix}, \tag{6}$$

satisfying the generating relations

$$R_1{}^2 = R_2{}^2 = R_3{}^2 = (R_1 R_3)^2 = (R_1 R_2)^3 = (R_2 R_3)^5 = 1. \tag{7}$$

Their representations in the $E_\perp$ space follows from (6) by algebraic conjugation. The orbits of the weights $v_i$ under the icosahedral group $H_3$ can be obtained by applying the group elements on the weight vectors as

$$\frac{(v_1)_{h_3}}{\sqrt{2}} = \frac{1}{2}\{(\pm 1, \pm\tau, 0), (\pm\tau, 0, \pm 1), (0, \pm 1, \pm\tau)\},$$

$$\frac{(\tau^{-1}v_2)_{h_3}}{2\sqrt{2}} = \frac{1}{2}\left\{(\pm 1,0,0),(0,0,\pm 1),(0,\pm 1,0),\frac{1}{2}(\pm 1,\pm\sigma,\pm\tau),\frac{1}{2}(\pm\sigma,\pm\tau,\pm 1),\frac{1}{2}(\pm\tau,\pm 1,\pm\sigma)\right\}, \quad (8)$$

$$\frac{(\tau^{-1}v_3)_{h_3}}{\sqrt{2}} = \frac{1}{2}\{(\pm 1,\pm 1,\pm 1),(0,\pm\tau,\pm\sigma),(\pm\tau,\pm\sigma,0),(\pm\sigma,0,\pm\tau)\},$$

where the notation $(v_i)_{h_3}$ is introduced for the set of vectors generated by the action of the icosahedral group on the weight $v_i$.

The sets of vectors in (8) represent the vertices of an icosahedron, an icosidodecahedron and a dodecahedron respectively in the $E_\parallel$ space. The weights $v_1$, $v_2$ and $v_3$ denote also the 5-fold, 2-fold and 3-fold symmetry axes of the icosahedral group. The union of the orbits of $\frac{(v_1)_{h_3}}{\sqrt{2}}$ and $\frac{(\tau^{-1}v_3)_{h_3}}{\sqrt{2}}$ constitute the vertices of a rhombic triacontahedron.

It is obvious that the Coxeter group $D_6$ is symmetric under the algebraic conjugation and this is more apparent in the characteristic equation of the Coxeter element of $D_6$ given by

$$(\lambda^3 + \sigma\lambda^2 + \sigma\lambda + 1)(\lambda^3 + \tau\lambda^2 + \tau\lambda + 1) = 0, \quad (9)$$

whose eigenvalues lead to the Coxeter exponents of $D_6$. The first bracket is the characteristic polynomial of the Coxeter element of the matrices in (6) describing it in the $E_\parallel$ space and the second bracket describes it in the $E_\perp$ space (Koca, Koc, Al-Barwani, 2001).

Therefore, projection of $D_6$ into either space is the violation of the algebraic conjugation. It would be interesting to discuss the projections of the fundamental polytopes of $D_6$ into 3D space possessing the icosahedral symmetry. It is beyond the scope of the present paper however we may discuss a few interesting cases.

The orbit generated by the weight $\omega_1$ is a polytope with 12 vertices called cross polytope and represents an icosahedron when projected into either space. The orbit of weight $\omega_2$ constitutes the "root polytope" of $D_6$ with 60 vertices which projects into 3D space as two icosidodecahedra with 30 vertices each, the ratio of radii of the circumspheres is $\tau$. The dual of the root polytope is the union of the three polytopes generated by $\omega_1$, $\omega_5$ and $\omega_6$ which constitute the Voronoi cell of the lattice $D_6$ as mentioned earlier and projects into an icosahedron and two rhombic triacontahedra. Actually, they consist of three icosahedra with the ratio of radii $1, \tau, \tau^2$ and two dodecahedra with the radii in proportion to $\tau$. The orbit generated by the weight vector $\omega_3$ is a polytope with 160 vertices and constitutes the Voronoi cell of the weight lattice $D_6^*$. It projects into two dodecahedra and two polyhedra with 60 vertices each. Voronoi cells can be used as windows for the cut and projects scheme however we prefer the direct projection of the root lattice as described in what follows.

A general root vector can be decomposed in terms of the weights $v_i(\acute{v}_i)$ as

$$m_1 l_1 + m_2 l_2 + m_3 l_3 + m_4 l_4 + m_5 l_5 + m_6 l_6 = \frac{1}{\sqrt{2+\tau}}[(m_1 - m_2 + \tau m_5 - \tau m_6)\, v_1 +$$
$$(m_2 - m_3 + \tau m_4 - \tau m_5)\, v_2 + (m_5 + m_6 + \tau m_3 - \tau m_4)\, v_3] + (v_i \rightarrow \acute{v}_i, \tau \rightarrow \sigma). \quad (10)$$

Projection of an arbitrary vector of $D_6$ into the space $E_\parallel$, or more thoroughly, onto a particular weight vector, for example, onto the weight $v_1$ is given by

$$[m_1 l_1 + m_2(l_2 + l_3 + l_4 + l_5 - l_6)]_\parallel \equiv [(m_1 - m_2)\,\omega_1 + 2m_2\,\omega_5]_\parallel = c\left(\frac{v_1}{\sqrt{2}}\right) \quad (11)$$

and represents an icosahedron where $c \equiv c(m_1, m_2) = \sqrt{\frac{2}{2+\tau}}(m_1 - m_2 + 2m_2\tau)$ is an overall scale factor. The subscript ∥ means the projection into the space $E_\parallel$. The expression in (11) shows that $m_1 + m_2 = even$ implying that the pair of integers $(m_1, m_2)$ are either both even or both odd. We will see that not only icosahedron but also dodecahedron and the icosahedral Archimedean polyhedra can be obtained by relations similar to (11). The Platonic and Archimedean polyhedra with icosahedral symmetry are listed in Table 1 as projections of $D_6$ lattice vectors determined by the pair of integers $(m_1, m_2)$. Table 1 shows that only a certain subset of vectors of $D_6$ project onto regular icosahedral polyhedra. We will see in the next section that the Danzer tiling also restricts the $D_6$ lattice in a domain where the vectors are determined by the pair of integers $(m_1, m_2)$. These quantities studied here are potentially applicable to many graph indices (Alhevaz, Baghipur, Shang, 2019)

**Table 1**
Platonic and Archimedean icosahedral polyhedra projected from $D_6$ (vertices are orbits under the icosahedral group).

| Vector of the Polyhedron in $H_3$ $c \equiv \sqrt{\frac{2}{2+\tau}}(m_1 - m_2 + 2m_2\tau)$ | | Corresponding Vector in $D_6$ $m_1, m_2 \in \mathbb{Z}$ |
|---|---|---|
| Icosahedron $c\frac{v_1}{\sqrt{2}} = \frac{c}{2}(1, \tau, 0)$ | 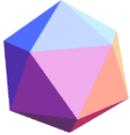 | $m_1 l_1 + m_2(l_2 + l_3 + l_4 + l_5 - l_6)$ $= (m_1 - m_2)\omega_1 + 2m_2\omega_5$ $m_1, m_2 \in 2\mathbb{Z}$ or $2\mathbb{Z} + 1$ |
| Dodecahedron $c\frac{v_3}{\sqrt{2}} = \frac{c}{2}(0, \tau^2, 1)$ | 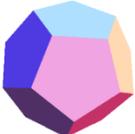 | $\frac{1}{2}[(m_1 + 3m_2)(l_1 + l_2 + l_3)$ $+ (m_1 - m_2)(l_4 + l_5 + l_6)]$ $= (m_1 - m_2)\omega_6 + 2m_2\omega_3$ $m_1, m_2 \in 2\mathbb{Z}$ or $2\mathbb{Z} + 1$ |
| Icosidodecahedron $c\frac{v_2}{\sqrt{2}} = c(0, \tau, 0)$ | 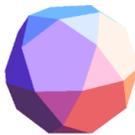 | $(m_1 + m_2)(l_1 + l_2) + 2m_2(l_3 + l_4)$ $= (m_1 - m_2)\omega_2 + 2m_2\omega_4$ $m_1, m_2 \in 2\mathbb{Z}$ or $2\mathbb{Z} + 1$ |
| Truncated Icosahedron $c\frac{(v_1 + v_2)}{\sqrt{2}} = \frac{c}{2}(1, 3\tau, 0)$ | 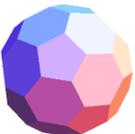 | $(2m_1 + m_2)l_1 + (m_1 + 2m_2)l_2$ $+ m_2(3l_3 + 3l_4 + l_5 - l_6)$ $= (m_1 - m_2)(\omega_1 + \omega_2) + 2m_2(\omega_4 + \omega_5)$ $m_1, m_2 \in 2\mathbb{Z}$ or $2\mathbb{Z} + 1$ |
| Small Rhombicosidodecahedron $c\frac{(v_1 + v_3)}{\sqrt{2}} = \frac{c}{2}(1, 2\tau + 1, 1)$ | 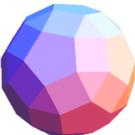 | $\frac{1}{2}[(3m_1 + 3m_2)l_1 + (m_1 + 5m_2)(l_2 + l_3)$ $+ (m_1 + m_2)(l_4 + l_5)$ $+ (m_1 - 3m_2)l_6]$ $= (m_1 - m_2)(\omega_1 + \omega_6) + 2m_2(\omega_3 + \omega_5)$ |
| Truncated Dodecahedron $$c\frac{(v_2+v_3)}{\sqrt{2}} = \frac{c}{2}(0, 3\tau+1, 1)$$ | 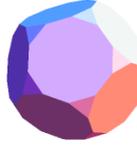 | $\frac{1}{2}[(3m_1+5m_2)(l_1+l_2) + (m_1+7m_2)l_3$ $+ (m_1+3m_2)l_4$ $+ (m_1-m_2)(l_5+l_6)]$ $= (m_1-m_2)(\omega_2+\omega_6) + 2m_2(\omega_3+\omega_4)$ |
| Great Rhombicosidodecahedron $$c\frac{(v_1+v_2+v_3)}{\sqrt{2}} = \frac{c}{2}(1, 4\tau+1, 1)$$ | 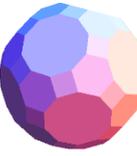 | $\frac{1}{2}[5(m_1+m_2)l_1 + (3m_1+7m_2)l_2 + (m_1+9m_2)l_3$ $+ (m_1+5m_2)l_4 + (m_1+m_2)l_5$ $+ (m_1-3m_2)l_6]$ $= (m_1-m_2)(\omega_1+\omega_2+\omega_6) + 2m_2(\omega_3+\omega_4+\omega_5)$ |

## 3. Danzer's *ABCK* tiles and $D_6$ lattice

We introduce the *ABCK* tiles with their coordinates in Fig. 2 as well as their images in lattice $D_6$. For a fixed pair of integers $(m_1, m_2) \neq (0,0)$ let us define the image of $K$ by vertices $D_1(m_1, m_2), D_2(m_1, m_2), D_3(m_1, m_2)$ and $D_0(0,0)$ in lattice $D_6$ and its projection in 3D space where $D_0(0,0)$ represents the origin. They are given as

$$D_1(m_1, m_2) = \frac{cv_1}{\sqrt{2}} = [m_1 l_1 + m_2(l_2 + l_3 + l_4 + l_5 - l_6)]_\parallel = [(m_1-m_2)\omega_1 + 2m_2\omega_5]_\parallel,$$

$$D_2(m_1, m_2) = \frac{cv_2}{2\sqrt{2}} = \frac{1}{2}[(m_1+m_2)(l_1+l_2) + 2m_2(l_3+l_4)]_\parallel = \frac{1}{2}[(m_1-m_2)\omega_2 + 2m_2\omega_4]_\parallel, \quad (12)$$

$$D_3(m_1, m_2) = \frac{c\tau^{-1}v_3}{\sqrt{2}} = \frac{1}{2}[(m_1+m_2)(l_1+l_2+l_3) + (-m_1+3m_2)(l_4+l_5+l_6)]_\parallel$$

$$= [2m_2\omega_3 + (-m_1+3m_2)\omega_6]_\parallel,$$

which follows from Table 1 with redefinitions of the pair of integers $(m_1, m_2)$. With removal of the notation $\parallel$ they represent the vectors in 6D space. A set of scaled vertices defined by $\acute{D}_i \equiv \frac{D_i(m_1,m_2)}{c}$ are used for the vertices of the Danzer's tetrahedra in Fig. 2.

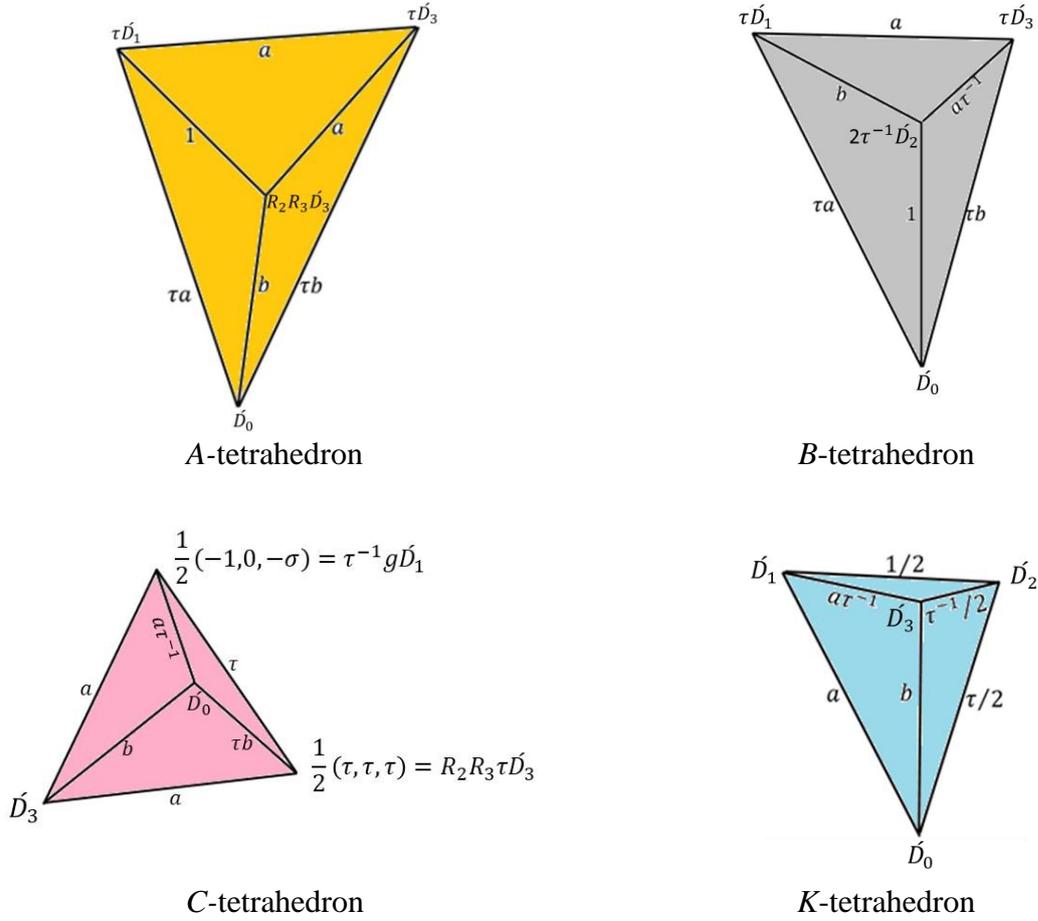

*A*-tetrahedron      *B*-tetrahedron

*C*-tetrahedron      *K*-tetrahedron

**Figure 2**
Danzer's tetrahedra.

It is instructive to give a brief introduction to Danzer's tetrahedra before we go into details. The edge lengths of $ABC$ tetrahedra are related to the weights of the icosahedral group $H_3$ whose edge lengths are given by $a = \frac{\sqrt{2+\tau}}{2} = \frac{\|v_1\|}{\sqrt{2}}$, $b = \frac{\sqrt{3}}{2} = \frac{\|\tau^{-1}v_3\|}{\sqrt{2}}$, $1 = \frac{\|\tau^{-1}v_2\|}{\sqrt{2}}$ and their multiples by $\tau$ and $\tau^{-1}$, where $a, b, 1$ are the original edge lengths introduced by Danzer. However, the tetrahedron $K$ has edge lengths also involving $\frac{1}{2}, \frac{\tau}{2}, \frac{\tau^{-1}}{2}$. As we noted in Sec. 2 the vertices of a rhombic triacontahedron consist of the union of the orbits $\frac{(v_1)_{h_3}}{\sqrt{2}}$ and $\frac{(\tau^{-1}v_3)_{h_3}}{\sqrt{2}}$. One of its cells is a pyramid based on a golden rhombus with vertices

$$\frac{v_1}{\sqrt{2}} = \frac{1}{2}(1, \tau, 0), \frac{\tau^{-1}v_3}{\sqrt{2}} = \frac{1}{2}(0, \tau, -\sigma), R_1 \frac{v_1}{\sqrt{2}} = \frac{1}{2}(-1, \tau, 0), R_3 \frac{\tau^{-1}v_3}{\sqrt{2}} = \frac{1}{2}(0, \tau, \sigma), \qquad (13)$$

and the apex is at the origin. The coordinate of the intersection of the diagonals of the rhombus is the vector $\frac{v_2}{2\sqrt{2}} = \frac{1}{2}(0, \tau, 0)$ and its magnitude is the in-radius of the rhombic triacontahedron. Therefore, the weights $\frac{v_1}{\sqrt{2}}, \frac{\tau^{-1}v_3}{\sqrt{2}}, \frac{v_2}{2\sqrt{2}}$ and the origin can be taken as the vertices of the tetrahedron $K$. As such, it is the fundamental region of the icosahedral group from which the rhombic triacontahedron is generated (Coxeter, 1973). The octahedra generated by these tetrahedra denoted by $<B>$, $<C>$ and $<A>$ comprise 4 copies of each obtained by a group of order 4 generated by two commuting generators $R_1$ and $R_3$ or their conjugate groups. The octahedron $<K>$ consists of $8K$ generated by a group of order 8 consisting of three commuting generators $R_1$, $R_3$ as mentioned earlier and the third generator $R_0$ is an affine reflection (Humphrey,1992) with respect to the golden rhombic face induced by the affine Coxeter group $\widetilde{D}_6$.

One can dissect the octahedron $<K>$ into three non-equivalent pyramids with rhombic bases by cutting along the lines orthogonal to three planes $v_1 - v_2, v_2 - v_3, v_3 - v_1$. One of the pyramids is generated by the tetrahedron $K$ upon the actions of the group generated by the reflections $R_1$ and $R_3$. If we call $R_1 K = \acute{K}$ as the mirror image of $K$ the others can be taken as the tetrahedra obtained from $K$ and $R_1 K$ by a rotation of $180^0$ around the axis $v_2$. A mirror image of $4K = (2K + 2\acute{K})$ with respect to the rhombic plane (corresponding to an affine reflection) complements it up to an octahedron of $8K$ as we mentioned above. In addition to the vertices in (13) the octahedron $8K$ also includes the vertices $(0,0,0)$ and $(0,\tau,0)$. Octahedral tiles are depicted in Fig. 3. Dissected pyramids constituting the octahedron $<K>$ as depicted in Fig. 4 have bases, two of which with bases of golden rhombuses with edge lengths (heights) $\tau^{-1}a\left(\frac{\tau}{2}\right), a\left(\frac{\tau^{-1}}{2}\right)$ and the third is the one with $b(\frac{1}{2})$. The rhombic triacontahedron consists of $60K + 60\acute{K} = 30(2K + 2\acute{K})$ where $(2K + 2\acute{K})$ form a cell of pyramid based on a golden rhombus of edge $\tau^{-1}a$ and height $\frac{\tau}{2}$. Faces of the rhombic triacontahedron are the rhombuses orthogonal to one of 30 vertices of the icosidodecahedron given in (8).

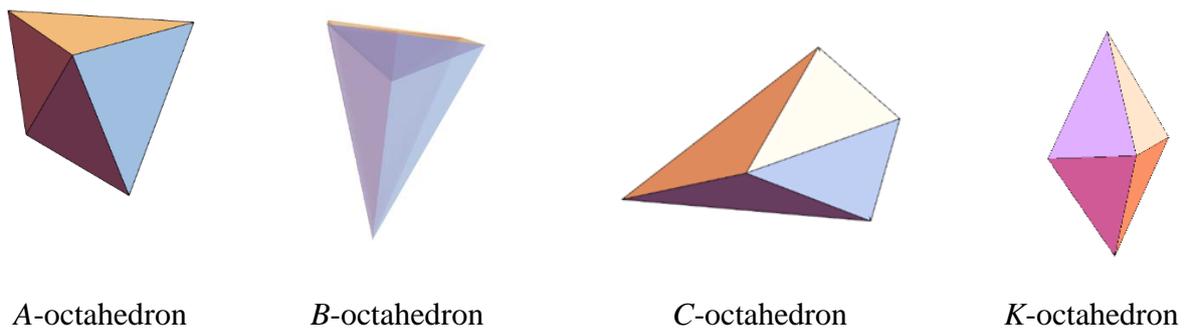

$A$-octahedron　　　$B$-octahedron　　　$C$-octahedron　　　$K$-octahedron

**Figure 3**
The octahedra generated by the $ABCK$ tetrahedra.

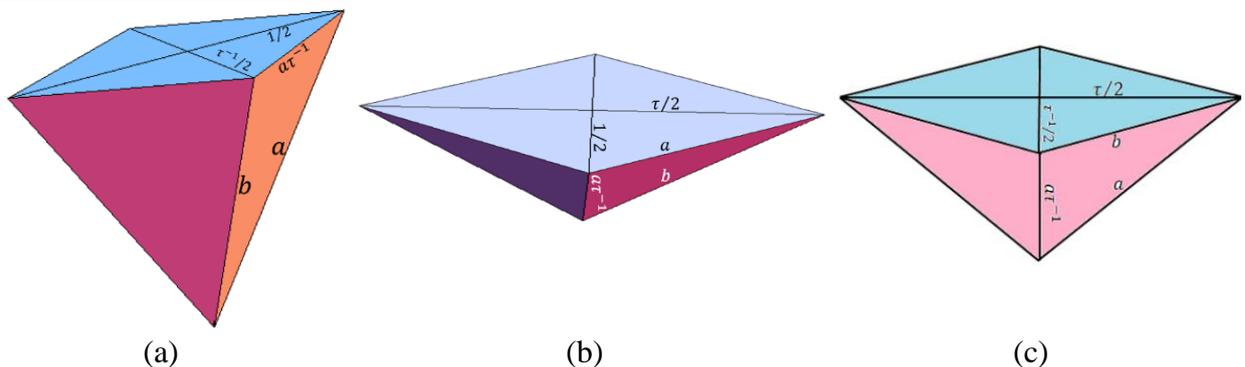

(a)　　　　　　　　　　(b)　　　　　　　　　　(c)

**Figure 4**
The three pyramids composed of $4K$ tetrahedra.

When the $4K$ of Fig. 4 (a) is rotated by the icosahedral group, it generates the rhombic triacontahedron as shown in Fig. 5.

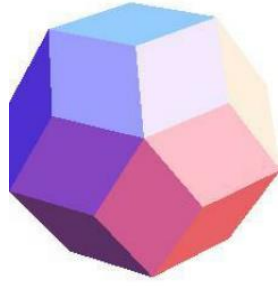

**Figure 5**
A view of rhombic triacontahedron.

The octahedron $<B>$ is a non-convex polyhedron whose vertices can be taken as $\tau \acute{D}_1$, $\tau R_1 \acute{D}_1, \tau \acute{D}_3, \tau R_3 \acute{D}_3, 2\tau^{-1}\acute{D}_2$ and the sixth vertex is the origin $\acute{D}_0$. It consists of 4 triangular faces with edges $a, \tau a, \tau b$ and 4 triangular faces with edges $\tau^{-1}a, a, b$. Full action of the icosahedral group on the tetrahedron $B$ would generate the "$B$-polyhedron" consisting of $60B + 60\acute{B}$ where $\acute{B} = R_1 B$ is the mirror image of $B$. It has 62 vertices (12 like $\frac{(\tau v_1)_{h_3}}{\sqrt{2}}$, 30 like $\frac{(\tau^{-1} v_2)_{h_3}}{\sqrt{2}}$, 20 like $\frac{(v_3)_{h_3}}{\sqrt{2}}$), 180 edges and 120 faces consisting of faces with triangles of edges with $\tau^{-1}a, a, b$. The face transitive "$B$-polyhedron" is depicted in Fig. 6 showing 3-fold and 5-fold axes simultaneously.

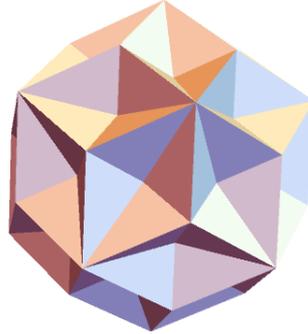

**Figure 6**
The "$B$-polyhedron".

The octahedron $<C>$ is a convex polyhedron which can be represented by the vertices

$$\tfrac{1}{2}(\tau, 0, 1), \tfrac{1}{2}(\tau, \tau, \tau), \tfrac{1}{2}(1, \tau, 0), \tfrac{1}{2}(\tau^2, 1, 0), \tfrac{1}{2}(\tau^2, \tau, 1), (0,0,0). \qquad (14)$$

This is obtained from that of Fig. 2 by a translation and inversion. It consists of 4 triangular faces with edges $a, b, \tau b$ and 4 triangular faces with edges $\tau^{-1}a, a$ and $b$.

The "$C$-polyhedron" is a non-convex polyhedron with 62 vertices (12 like $\frac{(v_1)_{h_3}}{\sqrt{2}}$, 30 like $\frac{(v_2)_{h_3}}{\sqrt{2}}$, 20 like $\frac{(v_3)_{h_3}}{\sqrt{2}}$), 180 edges and 120 faces consisting of only one type of triangular face with edge lengths $\tau^{-1}a, a$ and $b$ as can be seen in Fig. 7 (see also the reference Baake & Grimm, 2013, page 234). It is also face transitive same as "$B$-polyhedron".

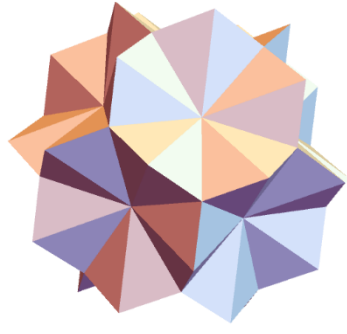

**Figure 7**
The "*C*-polyhedron".

The octahedron $<A>$ is a non-convex polyhedron which can be represented by the vertices

$$\tfrac{1}{2}(-\tau,0,1), \tfrac{1}{2}(1,-\tau,0), \tfrac{1}{2}(1,1,1), \tfrac{1}{2}(\sigma,0,-\tau), \tfrac{1}{2}(\sigma,-\tau,1), (0,0,0), \tag{15}$$

as shown in Fig. 3.

Now we discuss how the tiles are generated by inflation with an inflation factor $\tau$. First of all, let us recall that the projection of subset of the $D_6$ vectors specified by the pair of integers $(m_1, m_2)$ are some linear combinations of the weights $v_i, (i = 1, 2, 3)$ with an overall factor $c$. It is easy to find out the vector of $D_6$ corresponding to the inflated vertex of any *ABCK* tetrahedron by noting that $c(\tau^n v_i) = c(\acute{m}_1, \acute{m}_2) v_i$. Here we use $\tau^n = F_{n-1} + F_n \tau$ and we define

$$\acute{m}_1 \equiv m_1 F_{n-1} + \tfrac{1}{2}(m_1 + 5m_2)F_n, \quad \acute{m}_2 \equiv m_2 F_{n-1} + \tfrac{1}{2}(m_1 + m_2)F_n, \tag{16}$$

where $F_n$ represents the Fibonacci sequence satisfying

$$F_{n+1} = F_n + F_{n-1}, \quad F_{-n} = (-1)^{n+1} F_n, \quad F_0 = 0, F_1 = 1. \tag{17}$$

It follows from (16) that $\acute{m}_1 + \acute{m}_2 =$ even if $m_1 + m_2 =$ even, otherwise $\acute{m}_1, \acute{m}_2 \in \mathbb{Z}$. This proves that the pair of integers $(\acute{m}_1, \acute{m}_2)$ obtained by inflation of the vertices of the Danzer's tetrahedra remain in the subset of $D_6$ lattice. We conclude that the inflated vectors by $\tau^n$ in (12) can be obtained by replacing $m_1$ by $\acute{m}_1$ and $m_2$ by $\acute{m}_2$. For example radii of the icosahedra projected by $D_6$ vectors

$$2\omega_1, 2\omega_5, 2\omega_1 + 2\omega_5, 2\omega_1 + 4\omega_5 \tag{18}$$

are in proportion to $1, \tau, \tau^2$ and $\tau^3$ respectively. It will not be difficult to obtain the $D_6$ image of any general vector in the 3D space in the form of $c(\tau^p, \tau^q, \tau^r)$ where $p, q, r \in \mathbb{Z}$. Now we discuss the inflation of each tile one by one.

**Construction of $\tau K = B + K$**

The vertices of $\tau K$ is shown in Fig. 8 where the origin coincides with one of the vertices of *B* as shown in Fig. 2 and the other vertices of *K* and *B* are depicted. A transformation is needed to translate *K* to its new position. The face of *K* opposite to the vertex $\acute{D}_2$ having a normal vector $\frac{1}{2}(-1, -\sigma, -\tau)$ outward should match with the face of $B$ opposite to the vertex $\acute{D}_0$ with a normal vector $\frac{1}{2}(\sigma, \tau, -1)$ inward. For this reason, we perform a rotation of $K$ given by

$$g_K = \begin{bmatrix} 0 & -1 & 0 \\ 0 & 0 & -1 \\ 1 & 0 & 0 \end{bmatrix} \quad (19)$$

matching the normal of these two faces. A translation by the vector $\frac{1}{2}(\tau, \tau^2, 0)$ will locate $K$ in its proper place in $\tau K$. This is the simplest case where a rotation and a translation would do the work. The corresponding rotation and translation in $D_6$ can be calculated easily and the results are illustrated in Table 2.

**Table 2**
Rotation and Translation in $H_3$ and $D_6$.

| $H_3$ | | $D_6$ |
|---|---|---|
| $g_K$ | ← rotation→ | $g_K$: $(1\bar{3}6)(24\bar{5})$ |
| $t_K = \frac{1}{2}(\tau, \tau^2, 0)$ | ←translation→ | $\frac{1}{2}[(m_1 + 5m_2)l_1 + (m_1 + m_2)(l_2 + l_3 + l_4 + l_5 - l_6)]$ |

The rotation in $D_6$ is represented as the permutations of components of the vectors in the $l_i$ basis.

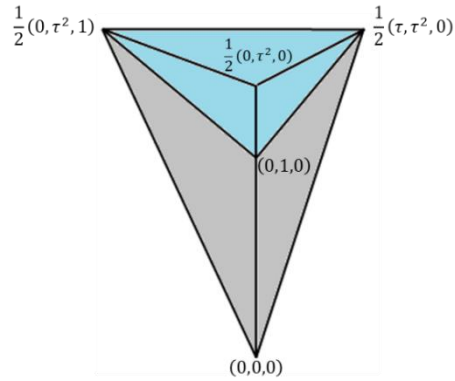

**Figure 8**
Vertices of $\tau K$.

Rhombic triacontahedron generated by $\tau K$ now consists of a "$B$-Polyhedron" centered at the origin and 30 pyramids of $4K$ with golden rhombic bases of edge lengths $a$ and heights $\frac{\tau^{-1}}{2}$ occupy the 30 inward gaps of "$B$-Polyhedron".

**Construction of $\tau B = C + 4K + B_1 + B_2$**

The first step is to rotate $B$ by a matrix

$$g_B = \frac{1}{2}\begin{bmatrix} -\sigma & \tau & -1 \\ -\tau & 1 & -\sigma \\ 1 & -\sigma & \tau \end{bmatrix}, \quad (20)$$

and then inflate by $\tau$ to obtain the vertices of $\tau g_B B$ as shown in Fig. 9.

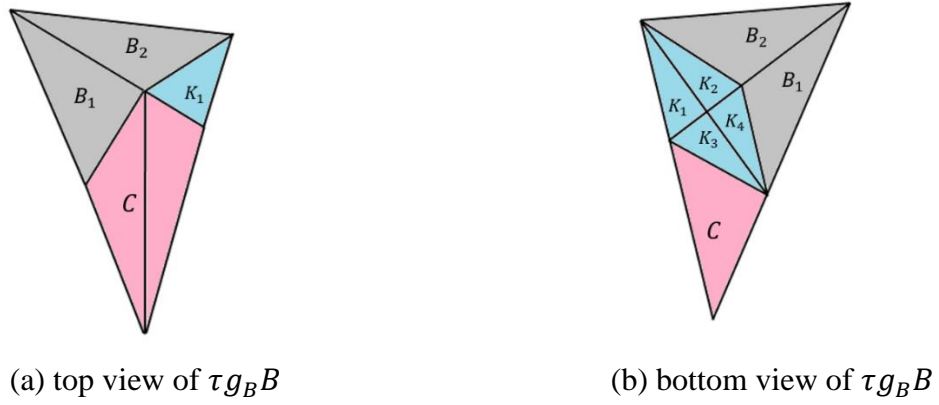

(a) top view of $\tau g_B B$          (b) bottom view of $\tau g_B B$

**Figure 9**
Views from $\tau g_B B$.

representation of $\tau g_B B$ is chosen to coincide the origin with the $\acute{D}_0$ vertex of tetrahedron C. It is obvious to see that the bottom view illustrates that the $4K$ is the pyramid based on the rhombus with edge length (height), $b(\frac{1}{2})$ as shown in Fig. 4 (c). The $4K$ takes its position in Fig. 9 by a rotation followed by a translation given by

$$\text{rotation: } g_{4K} = \frac{1}{2}\begin{bmatrix} 1 & -\sigma & -\tau \\ -\sigma & \tau & 1 \\ -\tau & 1 & \sigma \end{bmatrix}, \text{ translation: } t_{4K} = \frac{1}{2}(\tau, 0, 1). \qquad (21)$$

To translate $B_1$ and $B_2$ we follow the rotation and translation sequences as

$$\text{rotation: } g_{B_1} = \frac{1}{2}\begin{bmatrix} \sigma & -\tau & 1 \\ -\tau & 1 & -\sigma \\ -1 & \sigma & -\tau \end{bmatrix}, \text{ translation: } t_B = \frac{1}{2}(\tau^3, 0, \tau^2),$$

$$\text{rotation: } g_{B_2} = \frac{1}{2}\begin{bmatrix} -\sigma & -\tau & 1 \\ \tau & 1 & -\sigma \\ 1 & \sigma & -\tau \end{bmatrix}, \text{ translation: } t_B = \frac{1}{2}(\tau^3, 0, \tau^2). \qquad (22)$$

**Construction of $\tau C = K_1 + K_2 + C_1 + C_2 + A$**

We inflate by $\tau$ the tetrahedron $C$ with vertices shown in Fig. 10. Top view of $\tau C$ where one of the vertices of $K_1$ is at the origin is depicted in Fig.1. The vertices of the constituting tetrahedra are given as follows

$$K_1: \{(0,0,0), \tfrac{1}{2}(1,1,1), \tfrac{1}{2}(\tau, 0, 1), \tfrac{1}{4}(\tau^2, \tau, 1)\};$$

$$K_2: \{\tfrac{1}{2}(1,1,1), \tfrac{1}{2}(\tau, 0, 1), \tfrac{1}{4}(\tau^2, \tau, 1), \tfrac{1}{2}(\tau^2, \tau, 1)\};$$

$$C_1: \{\tfrac{1}{2}(1,1,1), \tfrac{1}{2}(\tau, 0, 1), \tfrac{1}{2}(\tau^2, \tau, 1), \tfrac{1}{2}(\tau^2, \tau^2, \tau^2)\}; \qquad (23)$$

$$C_2: \{\tfrac{1}{2}(\tau, 0, 1), \tfrac{1}{2}(\tau^2, \tau, 1), \tfrac{1}{2}(\tau^2, 0, \tau), \tfrac{1}{2}(\tau^2, \tau^2, \tau^2)\};$$

$$A: \{\tfrac{1}{2}(\tau^2,\tau,1), \tfrac{1}{2}(\tau^2,0,\tau), \tfrac{1}{2}(\tau^2,\tau^2,\tau^2), \tfrac{1}{2}(\tau^3,\tau^2,\tau)\}.$$

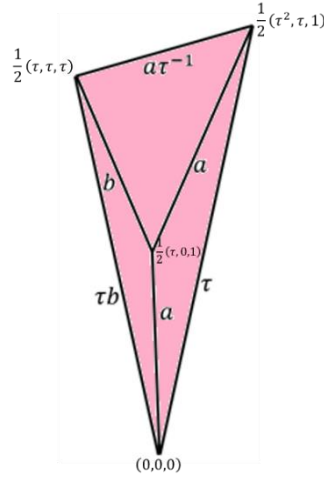

**Figure 10**
$C$-tetrahedron obtained from that of Fig. 2 by translation and inversion.

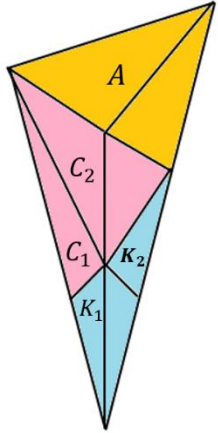 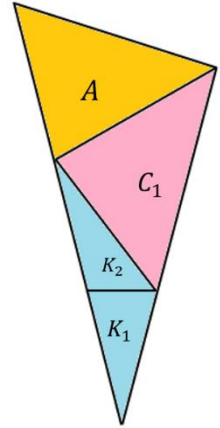

(a) Top view of $\tau C$          (b) Bottom view of $\tau C$

**Figure 11**
Views from $\tau C$.

The vertices of $K_1$ can be obtained from $K$ by a rotation $K_1 = g_{K_1} K$, where

$$g_{K_1} = \frac{1}{2}\begin{bmatrix} -\sigma & \tau & -1 \\ -\tau & 1 & -\sigma \\ 1 & -\sigma & \tau \end{bmatrix}, \quad (g_{K_1})^5 = 1. \tag{24}$$

The transformation to obtain $K_2$ is a rotation $g_{K_2}$ followed by a translation of $t_{K_2} = \tfrac{1}{2}(\tau^2,\tau,1)$ where

$$g_{K_2} = \frac{1}{2}\begin{bmatrix} -\sigma & -\tau & -1 \\ -\tau & -1 & -\sigma \\ 1 & \sigma & \tau \end{bmatrix}, \quad (g_{K_2})^{10} = 1. \tag{25}$$

To obtain the coordinates of $C_1$ and $C_2$ first translate $C$ by $-\frac{1}{2}(\tau,\tau,\tau)$ and rotate by $g_{C_1}$ and $g_{C_2}$ respectively followed by the translation $t_{C_1} = \frac{1}{2}(\tau^2, \tau^2, \tau^2)$ where

$$g_{C_1} = \frac{1}{2}\begin{bmatrix} \sigma & \tau & 1 \\ \tau & 1 & \sigma \\ 1 & \sigma & \tau \end{bmatrix}, \quad g_{C_2} = \frac{1}{2}\begin{bmatrix} 1 & -\sigma & -\tau \\ -\sigma & \tau & 1 \\ \tau & -1 & -\sigma \end{bmatrix}. \tag{26}$$

Similarly, vertices of tetrahedron $A$ is rotated by

$$g_A = \frac{1}{2}\begin{bmatrix} 1 & -\sigma & -\tau \\ -\sigma & \tau & 1 \\ -\tau & 1 & \sigma \end{bmatrix}, \tag{27}$$

followed by a translation $t_A = \frac{1}{2}(\tau^2, 0, \tau)$.

**Construction of $\tau A = 3B + 2C + 6K$**

It can also be written as $\tau A = C + K_1 + K_2 + B + \tau B$ where $\tau B$ is already studied. A top view of $\tau A$ is depicted in Fig.12 with the vertices of $\tau B$ given in Fig. 9.

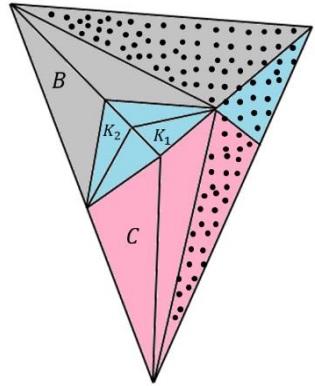

**Figure 12**
Top view of $\tau A$ (bottom view is the same as bottom view of $\tau B$- dotted region).

The vertices of $C$ $\{(0,0,0), \frac{1}{2}(\tau, 0,1), \frac{1}{2}(\tau^2, 1,0), \frac{1}{2}(\tau^2, \tau, 1)\}$ can be obtained from those in Fig. 2 by a translation $\frac{1}{2}(-\tau, -\tau, -\tau)$ followed by a rotation

$$g_C = \frac{1}{2}\begin{bmatrix} -1 & \sigma & -\tau \\ \sigma & -\tau & 1 \\ -\tau & 1 & -\sigma \end{bmatrix}. \tag{28}$$

Vertices of $K_1$ and $K_2$ are given by

$$K_1: \{\tfrac{1}{2}(\tau, 0,1), \tfrac{1}{2}(\tau^2, 1,0), \tfrac{1}{2}(\tau^2, \tau, 1), \tfrac{1}{4}(3\tau + 1, \tau, 1)\},$$

$$K_2: \{\tfrac{1}{2}(\tau,0,1), \tfrac{1}{2}(\tau^2,\tau,1), \tfrac{1}{4}(3\tau+1,\tau,1), \tfrac{1}{2}(2\tau,-\sigma,1)\}. \tag{29}$$

To obtain $K_1$ and $K_2$ with these vertices first rotate $K$ by $g_{K_1}$ and $g_{K_2}$ given by the matrices

$$g_{K_1} = \frac{1}{2}\begin{bmatrix} \sigma & \tau & -1 \\ \tau & 1 & -\sigma \\ 1 & \sigma & -\tau \end{bmatrix}, \quad g_{K_2} = \frac{1}{2}\begin{bmatrix} \sigma & \tau & 1 \\ \tau & 1 & \sigma \\ 1 & \sigma & \tau \end{bmatrix}, \tag{30}$$

and translate both by the vector $t_{K_1} = \tfrac{1}{2}(\tau,0,1)$. The vertices of $B$ in Fig.12 can be obtained by a rotation then a translation given respectively by

$$g_B = \frac{1}{2}\begin{bmatrix} -\tau & -1 & \sigma \\ -1 & -\sigma & \tau \\ -\sigma & -\tau & 1 \end{bmatrix}, \quad t_B = \tfrac{1}{2}(\tau^3, 0, \tau^2). \tag{31}$$

All these procedures are described in Table 3 which shows the rotations and translations both in $E_\parallel$ and $D_6$ spaces.

**Table 3**
Rotation and Translation in $H_3$ and $D_6$ (see the text for the definition of rotations in $E_\parallel$ space).

| $H_3$ | | | $D_6$ |
|---|---|---|---|
| | | $\tau K = B + K$ | |
| $g_K$ | ←rotation→ | | $(1\bar{3}6)(24\bar{5})$ |
| $t_K = \tfrac{c}{2}(\tau,\tau^2,0)$ | ←translation→ | | $\tfrac{1}{2}[(m_1+5m_2)l_1 + (m_1+m_2)(l_2+l_3+l_4+l_5-l_6)]$ |
| | | $\tau B = C + 4K + B_1 + B_2$ | |
| $g_B$ | ←rotation→ | | $(3)(126\bar{4}5)$ |
| $g_{4K}$ | ←rotation→ | | $(1)(23)(45)(6)$ |
| $t_{4K} = \tfrac{c}{2}(\tau,0,1)$ | ←translation→ | | $m_1 l_5 + m_2(l_1 - l_2 + l_3 - l_4 - l_6)$ |
| $g_{B_1}$ | ←rotation→ | | $(164 3\bar{5})(2)$ |
| $g_{B_2}$ | ←rotation→ | | $(1\bar{5}2\bar{1}\bar{5}2), (3\bar{6}4\bar{3}64)$ |
| $t_B = \tfrac{c}{2}(\tau^3,0,\tau^2)$ | ←translation→ | | $\tfrac{1}{2}[(3m_1+5m_2)l_5 + (m_1+3m_2)(l_1-l_2+l_3-l_4-l_6)]$ |
| | | $\tau C = K_1 + K_2 + C_1 + C_2 + A$ | |
| $g_{K_1}$ | ←rotation→ | | $(5\bar{4}621)(3)$ |
| | (identity translation) | | |
| $g_{K_2}$ | ←rotation→ | | $(1\bar{1})(6\bar{3}4\bar{5}2\bar{6}\bar{3}4\bar{5}2)$ |
| $t_K = \frac{c}{2}(\tau^2, \tau, 1)$ | ←translation→ | | $(m_1 + m_2)(l_1 + l_5)$ $+ 2m_2)(l_3 - l_6)$ |
| $t_C = -\frac{c}{2}(\tau, \tau, \tau)$ | ←translation→ | | $-\frac{1}{2}[(m_1 + 3m_2)(l_1 + l_3 + l_5)$ $+(m_1-m_2)(l_2 - l_4 - l_6)]$ |
| $g_{C_1}$ | ←rotation→ | | $(1)(4)(2\bar{6})(35)$ |
| $g_{C_2}$ | ←rotation→ | | $(1)(35\bar{6}42)$ |
| $t_C = \frac{c}{2}(\tau^2, \tau^2, \tau^2)$ | ←translation→ | | $(m_1 + 2m_2)(l_1 + l_3 + l_5)$ $+ m_2(l_2 - l_4 - l_6)$ |
| $g_A$ | ←rotation→ | | $(1)(6)(23)(45)$ |
| $t_A = \frac{c}{2}(\tau^2, 0, \tau)$ | ←translation→ | | $\frac{1}{2}[(m_1 + 5m_2)l_5$ $+(m_1 + m_2)(l_1 - l_2 + l_3 - l_4 - l_6)]$ |
| | | $\tau A = C + K_1$ $+ K_2 + B + \tau B$ | |
| $t_C = -\frac{c}{2}(\tau, \tau, \tau)$ | ←translation→ | | $-\frac{1}{2}[(m_1 + 3m_2)(l_1 + l_3 + l_5)$ $+(m_1 - m_2)(l_2 - l_4 - l_6)]$ |
| $g_C$ | ←rotation→ | | $(1\bar{1}\,(2\bar{4})(36)(5\bar{5})$ |
| $g_{K_1}$ | ←rotation→ | | $(1)\,(2543\bar{6})$ |
| $g_{K_2}$ | ←rotation→ | | $(1)(4)\,(2\bar{6})(35)$ |
| $t_K = \frac{c}{2}(\tau, 0, 1)$ | ←translation→ | | $m_1 l_5 +$ $m_2(l_1 - l_2 + l_3 - l_4 - l_6)$ |
| $g_B$ | ←rotation→ | | $(6\bar{5}\bar{1}6\bar{5}1), (\bar{4}2\bar{3}423)$ |
| $t_B = \frac{c}{2}(\tau^3, 0, \tau^2)$ | ←translation→ | | $\frac{1}{2}[(3m_1 + 5m_2)l_5 +$ $(m_1 + 3m_2)(l_1 - l_2 + l_3 - l_4 - l_6)]$ |

## 4. Discussions

The 6-dimensional reducible representation of the icosahedral subgroup of the point group of $D_6$ is decomposed into the direct sum of its two 3-dimensional representations described also by the direct sum of two graphs of icosahedral group. We have shown that the subset of the $D_6$ lattice characterized by a pair of integers $(m_1, m_2)$ with $m_1 + m_2 =$ even projects onto the Platonic and Archimedean polyhedra possessing icosahedral symmetry which are determined as the orbits of the fundamental weights $v_i$ or their multiples by $\tau^n$. The edge lengths of the Danzer tiles are related to the weights of the icosahedral group $H_3$ and via Table 1 to the weights of $D_6$, a group theoretical property which has not been discussed elsewhere. It turns out that the tetrahedron $K$ constitutes the fundamental region of the icosahedral group $H_3$ as being the cell of the rhombic triacontahedron. Images of the Danzer tiles and their inflations are determined in $D_6$ by employing translations and icosahedral rotations. This picture gives a one-to-one correspondence between the translations-rotations of 3D space and 6D space.

Faces of the Danzer tiles are all parallel to the faces of the rhombic triacontahedron; in other words, they are all orthogonal to the 2-fold axes. Since the faces of the Ammann rhombohedral tiles are orthogonal to the 2-fold axes of the icosahedral group it is this common feature that any tilings obtained from the Ammann tiles either with decoration or in the form of Socolar-Steinhardt model can also be obtained by the Danzer tiles. These geometrical properties have not been studied from the group theoretical point of view, a subject which, is beyond the present work but definitely deserves studying. Note that the inflation introduces a cyclic permutation among the rhombic triacontahedron, "$B$-polyhedron" and "$C$-polyhedron" as they occur alternatively centered at the origin.

Another novel feature of the paper is to show that *ABCK* tiles and their inflations are directly related to the transformations in the subset of $D_6$ lattice characterized by integers $(m_1, m_2)$ leading to an alternative projection technique different from the cut and project scheme. Details of the composition of the *ABCK* tiles into a structure with long-range quasiperiodic order follows from the inflation matrix (see for instance Baake & Grimm, 2013, pages 229-235).

**Acknowledgement**

We would like to thank Prof. Ramazan Koc for his contributions to Fig.6 and Fig.7.